\documentclass[notitlepage]{scrartcl}
\usepackage{enumitem}
\usepackage{amsmath}
\usepackage{amssymb}
\usepackage{amsthm}
\usepackage{abstract}
\usepackage{xcolor}
\usepackage{comment}
\usepackage{hyperref}
\hypersetup{
    colorlinks=true,
    linkcolor=blue,
    filecolor=magenta,      
    urlcolor=cyan
    }
\makeatletter
\newcommand*{\barfix}[2][.175ex]{%
  \mathpalette{\@barfix{#1}}{#2}%
}
\newcommand*{\@barfix}[3]{%
  \vbox{%
    \kern#1\relax
    \hbox{$#2#3\m@th$}%
  }%
}
\makeatother

\newcommand{\footremember}[2]{%
    \footnote{#2}
    \newcounter{#1}
    \setcounter{#1}{\value{footnote}}%
}

\usepackage[margin=1in]{geometry}
\title{\vspace{-2em}On the Performance of the Depth First Search Algorithm in Supercritical Random Graphs}
\author{%
Sahar Diskin \footremember{alley}{School of Mathematical Sciences, Tel Aviv University, Tel Aviv 6997801, Israel. Email: sahardiskin@mail.tau.ac.il.}%
\and Michael Krivelevich \footremember{trailer}{School of Mathematical Sciences, Tel Aviv University, Tel Aviv 6997801, Israel. Email: krivelev@tauex.tau.ac.il. Research supported in part by USA–Israel BSF grant 2018267 and by ISF grant
1261/17.}%
}
\begin{document}
	\maketitle
	
\begin{onecolabstract}
We consider the performance of the Depth First Search (DFS) algorithm on the random graph $G\left(n,\frac{1+\epsilon}{n}\right)$, $\epsilon>0$ a small constant. Recently, Enriquez, Faraud and M\'enard \cite{EFM} proved that the stack $U$ of the DFS follows a specific scaling limit, reaching the maximal height of $\left(1+o_{\epsilon}(1)\right)\epsilon^2n$. 
Here we provide a simple analysis for the typical length of a maximum path discovered by the DFS.
\end{onecolabstract}

\section{Introduction}
We consider the structure of the spanning tree of the giant component of $G(n,p)$ uncovered by the Depth First Search (DFS) algorithm, for the supercritical regime $p=\frac{1+\epsilon}{n}.$

As for the notation of the sets in the DFS algorithm, we follow the conventions similar to \cite{K and S}: 
We denote by $S$ the set of vertices whose exploration is complete; by $T$ the set of vertices not yet visited, and by $U$ the set of vertices which are currently being explored, kept in a stack. At any moment $0\le m\le {n \choose 2}$ in the DFS, we denote by $S(m)$, $T(m)$ and $U(m)$ the sets $S$, $T$ and $U$ (respectively) at $m.$

The algorithm starts with $S=U=\emptyset$ and $T=V(G)$, and ends when $U\cup T=\emptyset$. At each step, if $U$ is nonempty, the algorithm queries $T$ for neighbours of the last vertex in $U.$ The algorithm is fed $X_i$, $0\le i\le {n\choose 2}$, i.i.d Bernoulli$(p)$ random variables, each corresponding to a positive (with probability $p$) or negative (with probability $1-p$) answer to such a query. If $U$ is nonempty and the last vertex in $U$ has no more queries to ask, then we move the last vertex of $U$ to S. If $U=\emptyset$, we move the next vertex from $T$ into $U.$ Formally, after completing the discovery of all the connected components, we query all the remaining pairs of vertices that have not been queried by the DFS.

Enriquez, Faraud and M\'enard provided in \cite{EFM} an analysis of the performance of DFS: tracking the stack $U,$ they showed it follows a specific scaling limit, reaching the maximal height of $\left(1+o_{\epsilon}(1)\right)\epsilon^2n.$ Here we provide a simpler, and perhaps more telling argument for the typical maximal length of a path found by DFS. 

Our result is as follows:
\paragraph{Theorem 1} \textit{Let $\epsilon>0$ be a small enough constant, and let $p=\frac{1+\epsilon}{n}.$ Run the DFS algorithm on $G(n,p)$. Then, \textbf{whp}, a longest path in the obtained spanning forest is of length $\epsilon^2n+O(\epsilon^3)n.$\\}

We should note that while the precise length of a longest path in $G(n,p)$ is an open problem, it is known that a longest path is \textbf{whp} at least of length $\frac{4\epsilon^2}{3}n$ and at most $\frac{7\epsilon^2}{4}n$ (see \cite{Kemkes and Wormald}, \cite{Lucak}). Hence, while the DFS finds a path of the correct magnitude ($\Theta(\epsilon^2)n$) as was shown already in \cite{K and S}, the longest path found by the algorithm is significantly shorter than a longest path in the graph.

Furthermore, while we treat $\epsilon$ as a constant, our statements and proof hold for any $\epsilon=\epsilon(n)$ that tends to $0$ with $n\to\infty$, as long as $\epsilon(n)\gg n^{-1/3+o(1)}$ (see the comment following the proof of Lemma 2.3), covering a substantial part of the barely-supercritical regime as well.

\section{Two-step Analysis}
We define the excess of a connected graph $G=(V,E)$ to be $|E(G)|-|V(G)|+1.$ We define the excess of a graph to be the sum of the excesses of its connected components. 

We require the following well-known facts regarding $G(n,p)$ (see, for example, \cite{Intro to Random Graphs}):
\paragraph{Theorem 2.1} \textit{Let $\epsilon>0$ be a small enough constant. Then, \textbf{whp}: 
\begin{description}
    \item[1.] In $G\left(n, \frac{1+\epsilon}{n}\right)$ there is a unique giant component, $L_1$, whose size is asymptotic to $\Theta(\epsilon)n.$ All the other components are of size $O\left(\ln n/\epsilon^2\right)$.
    \item[2.] The excess of $G\left(n,\frac{1+\epsilon}{n}\right)$ is at most $6\epsilon^3n.$
    \item[3.] In $G\left(n, \frac{1-\epsilon}{n}\right)$, all the components are of size $O\left(\ln n/\epsilon ^2\right)$.
\end{description}}
When $p=\frac{1+\epsilon}{n}$, we call $G(n,p)$ a \textit{supercritical random graph}. When $p=\frac{1-\epsilon}{n}$ we call $G(n,p)$ a \textit{subcritical random graph}.

We also require the following simple lemma:
\paragraph{Lemma 2.2} \textit{Let $\epsilon>0$ be a small enough constant, and let $p=\frac{1+\epsilon}{n}.$ Then, \textbf{whp}, by the moment $m=n\ln^2 n$ we are already in the midst of discovering the giant component.}
\begin{proof}
By Theorem 2.1, the largest component is \textbf{whp} of size  $\Theta(\epsilon)n$, and all the other components are of size $O\left(\frac{\ln n}{\epsilon^2}\right)$. As long as we are prior to the discovery of the giant component, every time $U$ empties, the new vertex about to enter $U$ has probability at least $\Theta(\epsilon)$ to belong to the giant component. Every time a vertex that does not belong to the giant enters $U,$ $U$ empties after at most $O\left(n\frac{\ln n}{\epsilon^2}\right)$ queries, corresponding to at most $O\left(\frac{\ln n}{\epsilon^2}\right)$ positive answers. Therefore, the probability that after $n\ln^2 n$ rounds we are still not in the midst of discovering the giant component is at most $\left(1-\Theta(\epsilon)\right)^{\frac{n\ln^2n}{O\left(n\frac{\ln n}{\epsilon^2}\right)}}=\left(1-\Theta(\epsilon)\right)^{\Omega(\epsilon^2)\ln n}=o(1)$.
\end{proof}

We will focus on the stack of the DFS, $U,$ and its development throughout the DFS run. 
\subsection{The Straightforward Analysis}
In hindsight, we know that $U$ reaches its maximal height around the moment $\frac{\epsilon n^2}{1+\epsilon}.$ However, around this moment issues with critically begin to occur. We thus define two moments which will be useful as points of reference for us:
\begin{align}
    m_1:=\frac{\left(\epsilon-\epsilon^2\right)n^2}{1+\epsilon}, \qquad m_2:=\frac{(\epsilon-\epsilon^2+\epsilon^3)n^2}{1+\epsilon}.
\end{align}

The following straightforward lemma gives a bound on the height of $U$ at the moment $m_1$, depending only on the number of queries between $U$ and $T$, which we will analyse afterwards:
\paragraph{Lemma 2.3} \textit{Let $\epsilon>0$ be a small enough constant and let $p=\frac{1+\epsilon}{n}.$ Let $m_1$ be as defined in $(1)$. Run the DFS algorithm on $G(n,p)$. Then, at the moment $m_1$ we have \textbf{whp}:
$$|U(m_1)|=\frac{\epsilon^2n}{2}+\frac{q_{m_1}(U,T)}{n}+O(\epsilon^3)n,$$
where $q_{m_1}(U,T)$ is the number of queries between the vertices of $U(m_1)$ and $T(m_1)$ by moment $m_1.$
}
\begin{proof}
We consider the different types of queries that occurred by moment $m_1$:
\begin{description}
    \item[1.] $q_{m_1}(S,T)$ is the number of queries between the vertices in $S(m_1)$ and $T(m_1)$ by the moment $m_1.$ By properties of the DFS, 
    $$q_{m_1}(S,T)=|S(m_1)||T(m_1)|.$$
    \item[2.] $q_{m_1}(S\cup U)$ is the number of queries inside $S(m_1)\cup U(m_1)$ by the moment $m_1.$ By Theorem 2.1, the excess of the graph is \textbf{whp} at most $6\epsilon^3n.$ Hence, we have that \textbf{whp}:
    \begin{align*}
        {{|S(m_1)|+|U(m_1)|}\choose 2}-6\epsilon^3n^2\le q_{m_1}(S\cup U)\le {{|S(m_1)|+|U(m_1)|}\choose 2}.
    \end{align*}
    Indeed, there are ${{|S(m_1)|+|U(m_1)|}\choose 2}$ possible queries inside $S(m_1)\cup U(m_1)$. In order to obtain the full description of the graph, we will need to ask all these queries. Should there be more than $6\epsilon^3n^2$ queries remaining after the DFS run, there would be \textbf{whp} (by a standard Chernoff-type bound, see, for example, Theorem A.1.11 of \cite{Chernoff from K and S}) at least $6\epsilon^3n$ additional edges, contradicting Theorem 2.1.
    \item[3.] $q_{m_1}(U,T)$ is the number of queries between the vertices in $U(m_1)$ and $T(m_1)$ by the moment $m_1.$ 
\end{description}
These types of queries account for all the queries by moment $m_1.$ We thus have that:
\begin{align*}
    m_1=\frac{(\epsilon-\epsilon^2)n^2}{1+\epsilon}=q_{m_1}(S,T)+q_{m_1}(S\cup U)+q_{m_1}(U,T),
\end{align*}
and
\begin{align*}
    \Bigg|\left(q_{m_1}(S,T)+q_{m_1}(S\cup U)\right)-\left(|T(m_1)||S(m_1)|+\frac{\left(|S(m_1)|+|U(m_1)|\right)^2}{2}\right)\Bigg|\le 6\epsilon^3n^2.
\end{align*}
By Lemma 2.2, by the moment $n\ln^2n$ we are already in the midst of discovering the largest component. As such, by the moment $m_1$, $U$ emptied \textbf{whp} at most $2\ln^2n$ times (every time $U$ emptied we must have had at least $(1-\Theta(\epsilon))n$ queries, \textbf{whp}). Therefore, by properties of the DFS run and by Lemma 2.2 we have that \textbf{whp},
\begin{align*}
    \Bigg||S(m_1)|+|U(m_1)|-\sum_{i=1}^{m_1} X_i\Bigg|\le 2\ln^2n,
\end{align*}
and $|T(m_1)|=n-|S(m_1)|-|U(m_1)|.$  Using a standard Chernoff-type bound together with the union bound, we obtain that with exponentially high probability:
\begin{align*}
    \Bigg|\sum_{i=1}^{m_1}X_i-(\epsilon-\epsilon^2)n\Bigg|\le \epsilon^3n.
\end{align*}
Hence \textbf{whp},
\begin{align*}
    |S(m_1)|+|U(m_1)|=(\epsilon-\epsilon^2)n+O(\epsilon^3)n,
\end{align*}
and thus \textbf{whp},
\begin{align*}
    \frac{(\epsilon-\epsilon^2)n^2}{1+\epsilon}&=|T(m_1)||S(m_1)|+\binom{|S(m_1)+|U(m_1)|}{2}+q_{m_1}(U,T)+O(\epsilon^3)n\\
    &= (n-(\epsilon-\epsilon^2)n)\left((\epsilon-\epsilon^2) n-|U(m_1)|\right)+\frac{\epsilon^2 n^2}{2}+q_{m_1}(U,T)+O(\epsilon^3)n^2\\
    &= \epsilon n^2-\frac{3\epsilon^2n^2}{2}-n|U(m_1)|+q_{m_1}(U,T)+O(\epsilon^3)n^2,
\end{align*}
where the last equality follows since $U(m_1)$ spans a path, and \textbf{whp} a longest path is of length at most $2\epsilon^2n$ (see \cite{Lucak}).
Multiplying both sides of the inequality by $\frac{1+\epsilon}{n}$, we obtain that \textbf{whp}:
\begin{align*}
    \epsilon n-\epsilon^2n&= (1+\epsilon)\left(\epsilon n-\frac{3\epsilon^2 n}{2}-|U(m_1)|+\frac{q_{m_1}(U,T)}{n}+O(\epsilon^3)n\right)\\
    &=\epsilon n-\frac{\epsilon ^2 n}{2}-|U(m_1)|+\frac{q_{m_1}(U,T)}{n}+O(\epsilon^3)n,
\end{align*}
for small enough $\epsilon$. Rearranging, we derive that \textbf{whp}:
\begin{align*}
    |U(m_1)|=\frac{\epsilon^2 n}{2}+\frac{q_{m_1}(U,T)}{n}+O(\epsilon^3)n,
\end{align*}
as required.
\end{proof}
We remark that with slight adjustment in the proof of Lemma 2.2, we have that \textbf{whp} by the moment $\frac{n\ln^2n}{\epsilon}$ we are already in the midst of discovering the largest component. Then, with a more careful treatment of the error terms, the proof of Lemma 2.3 follows through for any $\epsilon\gg n^{-1/3+o(1)}$ (and subsequently, so do the proofs of the following lemmas and Theorem 1). 

An immediate corollary of Lemma 2.3 is that the DFS uncovers \textbf{whp} a path of size at least $\frac{\epsilon^2n}{2}-O(\epsilon^3)n.$ In order to obtain tight bounds, we will need to analyse the quantity $q_{m_1}(U,T)$.

\subsection{Estimating $q_{m_1}(U,T)$}

We now want to obtain a good estimate for $q_{m_1}(U,T)$. For that, we first observe that $G[T(m)]$ behaves like a random graph. Specifically, for $m\le m_1$, $G[T(m)]$ behaves like a supercritical random graph, having a unique giant component with all other components of size at most logarithmic in $n$; for $m\ge m_2$, $G[T(m)]$ behaves like a subcritical random graph, with all components of size at most logarithmic in $n$. For $m_1< m< m_2$, $G[T(m)]$ might behave like a critical random graph, however, these two moments are close enough so this does not affect the size of $U$ significantly. We now state and prove this formally:

\paragraph{Lemma 2.4} \textit{Let $\epsilon>0$ be a small enough constant. Let $p=\frac{1+\epsilon}{n}$, and let $m_1, m_2$ be as defined in $(1)$. Run the DFS on $G(n,p)$. Then, \textbf{whp}, for all $m\le m_1$, $G[T(m)]$ behaves like a supercritical random graph, and for all $m\ge m_2$, $G[T(m)]$ behaves like a subcritical random graph.}
\begin{proof}
First we note that since at any moment $m$ the vertices in $T(m)$ have not been queried against each other, $G[T(m)]$ is distributed like $G\left(|T(m)|, \frac{1+\epsilon}{n}\right)$ random graph. Now, let $f(\epsilon), g(\epsilon)$ be positive constants depending on $\epsilon.$ Then, $G[T(m)]$ is supercritical if $|T(m)|p\ge 1+f(\epsilon)$, and subcritical if $|T(m)|p\le 1-g(\epsilon)$. Recall that $|T(m)|=n-|S(m)|-|U(m)|$, and that by Lemma 2.2 and by a Chernoff-type bound, \textbf{whp} $$\Bigg||S(m)+|U(m)|-\sum_{i=1}^{m}X_i\Bigg|\le \ln^2n.$$
Substituting $m=m_1$, we have \textbf{whp} that:
    \begin{align*}
        |T(m_1)|p&\ge \left(n-(\epsilon-\epsilon^2)n-4\sqrt{n\ln n}\right)\frac{1+\epsilon}{n}\\
        &\ge 1+\epsilon^3-5\sqrt{\frac{\ln n}{n}}.
    \end{align*}
Similarly, substituting $m=m_2$ we get \textbf{whp} that
    \begin{align*}
        |T(m_2)|p&\le \left(n-(\epsilon-\epsilon^2+\epsilon^3)n+3\sqrt{n\ln n}\right)\frac{1+\epsilon}{n}\\
        &\le 1-\epsilon^4+4\sqrt{\frac{\ln n}{n}}.
    \end{align*}
All that is left is to note that, by properties of the DFS, for any two moments $m\le m'$ we have that $T(m')\subseteq T(m)$ and thus $|T(m')|\le |T(m)|$.
\end{proof}

We are now ready to provide a good estimate for $q_{m_1}(U,T)$:

\paragraph{Lemma 2.5} \textit{Let $\epsilon>0$ be a small enough constant. Let $p=\frac{1+\epsilon}{n}$, and let $m_1$ be as defined in $(1)$. Run the DFS on $G(n,p)$. Then, \textbf{whp},
\begin{align*}
    \frac{|U(m_1)|}{2}-8\epsilon^3n\le \frac{q_{m_1}(U,T)}{n} \le \frac{(1+\epsilon)|U(m_1)|}{2}.
\end{align*}}
\begin{proof}
At any moment $m\le m_1$,  by Lemma 2.4 $G[T(m)]$ behaves like a supercritical random graph. As such, by Theorem 2.1, \textbf{whp} it has a unique giant component of size linear in $n$, with all other components of size at most logarithmic in $n$. 

Consider a vertex that entered $U$ at some moment $m \le m_1$. If it belonged to the giant component of $G[T(m-1)]$, then we will explore all of the giant component of $G[T(m-1)]$, whose size is linear in $n$, before it will move out of $U$. If it did not belong to the giant component of $G[T(m-1)]$, then we will explore a component of size logarithmic in $n$, before removing it from $U$. As such, all but the last $\ln^2 n$ vertices of $U(m_1)$ entered $U$ from a giant component (indeed, the last $\ln^2n$ vertices of $U(m_1)$ form a path, and a path of length $\ln^2n$ belongs to the giant component), and we can focus on these vertices.

Consider such a moment $m\le m_1$ where a vertex belonging to the giant component of $G[T(m-1)]$ entered $U$, and denote the last vertex in $U(m)$ by $v.$ Noting that these giant components are nested, and since by Lemma 2.4 \textbf{whp} $G[T(m_1)]$ has a giant component, we have that \textbf{whp} this holds for all $m\le m_1$. Hence, \textbf{whp} $G[T(m)]$ also has a giant component, and since $v$ belonged to the giant component of $G[T(m-1)]$, it must have at least one neighbour in the giant component of $G[T(m)].$ Let $q(v, m)$ be the random variable representing the number of queries the vertex $v$ in $U$ had against the vertices in $T(m)$, before the next vertex belonging to the giant of $G[T(m)]$ enters $U.$ 

For the upper bound, observe that $q(v, m)$ is stochastically dominated by the random variable $Uni(1, n)$, since we know that there is at least one neighbour of $v$ in the giant of $G[T(m)]$, and there are at most $n$ vertices in $T(m)$. Therefore, $q_{m_1}(U,T)$ is stochastically dominated by the sum of $|U(m_1)|$ $i.i.d$ random variables distributed according to $Uni(1,n)$, together with at most $n\ln^2 n$ additional queries accounting for the last $\ln^2n$ vertices in $U(m_1)$. By the Law of Large Numbers, we have that: 
$$P\left[\frac{q_{m_1}(U,T)}{n}\ge\frac{(1+\epsilon)|U(m_1)|}{2}\right]=o(1),$$
since $|U(m_1)|\ge \frac{\epsilon^2n}{2}.$

For the lower bound, observe that any additional neighbours that $v$ may have in the giant component of $G[T(m)]$, besides the one guaranteed by construction, contribute to the excess of the giant component. Indeed, the edges between $v$ and these additional neighbours will not be queried during the DFS run, since the entire giant component of $G[T(m)]$ will be explored before we return to $v$ in $U.$ By Theorem 2.1, the excess of the giant component is \textbf{whp} at most $6\epsilon^3n.$ Furthermore, while it is possible that some vertices moved from $T$ to $U$ (and later on to $S$) between the moment $m$ and the moment where we found the first neighbour in the giant, we still have that for all $m\le m_1$ \textbf{whp} $|T(m)|\ge |T(m_1)|\ge (1-2\epsilon)n.$ Thus $q_{m_1}(U,T)$ stochastically dominates the sum of $|U(m_1)|-6\epsilon^3n-\ln^2 n$ random variables distributed according to $Uni(1,(1-2\epsilon)n)$. Since $|U(m_1)|\ge \frac{\epsilon^2n}{2}$, by the Law of Large numbers we obtain the required lower bound \textbf{whp}. 
\end{proof}

\section{Proof of Theorem 1}
By Lemma 2.3 and Lemma 2.5, \textbf{whp} at the moment $m_1$ as defined in $(1)$,
\begin{align*}
    |U(m_1)|&=\frac{\epsilon^2n}{2}+\frac{q_{m_1}(U,T)}{n}+O(\epsilon^3)n\\
    &=\frac{\epsilon^2n}{2}+\frac{|U(m_1)|}{2}+O(\epsilon^3)n.
\end{align*}
Rearranging, we obtain that \textbf{whp} $|U(m_1)|=\epsilon^2n+O(\epsilon^3)n.$ 
This immediately proves the lower bound. For the upper bound, observe that by Lemma 2.4, between $m_1$ and $m_2$ (as defined in $(1)$) we have at most $O(\epsilon^3)n^2$ queries, corresponding to at most $O(\epsilon^3) n$ additional vertices to $U,$ \textbf{whp}. Afterwards, by Lemma 2.4, \textbf{whp} the DFS enters the subcritical phase, and by Theorem 2.1 \textbf{whp} all the components in $G[T]$ are of size logarithmic in $n$, at most. As such, $|U|$ could increase by at most $\ln^2n$, before decreasing back again. \qedsymbol

\end{document}